\documentclass{amsart}
\usepackage{amssymb,latexsym,amsfonts,amsthm,amsmath}
\usepackage{enumerate, epsfig}

\numberwithin{equation}{section}

\theoremstyle{plain}
\newtheorem{theorem}{Theorem}

\begin{document}
\title[Riemann surfaces with large corona constants]{An explicit example of Riemann surfaces with large bounds on corona solutions}
%\title{An explicit example of Riemann surfaces with large bounds on corona solutions}
\author{Byung-Geun Oh}
\address{Department of Mathematics, University of Washington, Box 354350, Seattle, WA 98195-4350}
\email{bgoh@math.washington.edu}
\date{June 12, 2005}
%\maketitle

\begin{abstract}
By modifying Cole's example, we construct explicit Riemann surfaces with large bounds on corona
solutions in an elementary way. 
\end{abstract}

\maketitle

\section{Introduction}
For a given Riemann surface $R$, we consider the algebra $H^\infty (R)$ of 
bounded analytic functions on $R$ separating the points in $R$. The corona problem 
asks whether $\iota (R)$ is dense in the maximal ideal space $\mathcal{M} (R)$ of $H^\infty (R)$,
where $\iota : R \to \mathcal{M} (R)$ is the natural inclusion defined by the point evaluation. 
If $\iota (R)$ is dense in $\mathcal{M} (R)$, we say that the \emph{corona theorem holds} for $R$. 
Otherwise $R$ is said to \emph{have corona}.

The corona theorem holds for $R$ if and only if the following statement is true (cf. Chap.~4 of \cite{Gam} or Chap.~VIII of \cite{Gar}):
for every collection $F_1, \ldots, F_n \in H^\infty (R)$ and any $\delta \in (0,1)$ with 
the property
\begin{equation}\label{bound}
  \delta \leq \max_j |F_j(\zeta)| \leq 1 \quad \mbox{for all } \zeta \in R,
\end{equation}
there exist $G_1, \ldots, G_n \in  H^\infty (R)$ such that
\begin{equation}\label{ideal}
 F_1 G_1 + F_2 G_2 + \cdots + F_n G_n =1.
\end{equation}

We refer to $G_1, \ldots, G_n$ as \emph{corona solutions}, $F_1, \ldots, F_n$ as 
\emph{corona data}, and $\max \{ \|G_1\|, \ldots, \|G_n\| \}$ 
as \emph{bound on the corona solutions} or \emph{corona constant}.
Here the notation $\| \cdot \|$ indicates the uniform norm. Throughout this paper, 
we assume that the corona data satisfies \eqref{bound} for the given $\delta$, where
the letter $\delta$ is reserved only for the above use. 

\begin{theorem}[B. Cole, cf. \cite{Gam}]\label{T}
For any $\delta \in (0,1)$ and $M >0$, there exist a finite bordered Riemann surface $R$
and corona data $F_1, F_2 \in H^\infty (R)$ such that any corona solutions $G_1, G_2 \in H^\infty (R)$
have a bound at least $M$; i.e., $\max \{ \|G_1\|, \|G_2\| \} \geq M$.
\end{theorem}

The purpose of this paper is to construct the Riemann surface $R$ in Theorem~\ref{T} in an
elementary way and describe it explicitly. Once Theorem~\ref{T} is proved,
it is possible to construct a Riemann surface with corona.

\begin{theorem}[B. Cole, cf. Theorem~4.2 of \cite{Gam}, p.~49--52]\label{TT}
There exists an open Riemann surface with corona.
\end{theorem}

The basic idea of the proof of Theorem~\ref{TT} 
is the fact that if a Riemann surface $R$ is obtained by
connecting two Riemann surfaces $R_1$ and $R_2$ with a thin strip, then any holomorphic function
on $R$ behaves almost independently on $R_1$ and $R_2$, respectively. 

The corona theorem holds for the unit disc \cite{Ca}, finitely connected domains in
$\mathbb{C}$ \cite{Gam2}, Denjoy domains \cite{GJ}, and various other classes of planar domains and
Riemann surfaces \cite{Al}, \cite{Be1}, \cite{Be2}, \cite{JM}, \cite{St}. 
On the other hand, examples of Riemann surfaces with corona other than Cole's example can be found 
in \cite{BD} and \cite{Ha}. Furthermore, by modifying the proof of Theorem~\ref{TT},
Cole's example can be used to obtain a Riemann surface with corona that is 
of Parreau-Widom type\footnote{A Riemann surface $R$ is called Parreau-Widom type if 
$\sum_{z \in E} G(z, w) < + \infty$, where $G(\cdot,w)$ is the Green's function on $R$ with the pole $w$
and $E = \{ z : \nabla G(z,w) = 0 \}$.} \cite{Na}.

The corona problem for a general domain 
in $\mathbb{C}$ is still open, and the answer is also unknown 
for a polydisc or a unit ball in $\mathbb{C}^n$, $n \geq 2$.

\section{Proof of Theorem~\ref{T}}\label{proof}
For given $\delta \in (0,1)$ and $M > 0$, we choose a natural number $n$ so that 
$\delta^n \leq \min \{ (16M)^{-1}, 1/4 \}$.
Let $d = 4 \delta^{n^2 + n}$ and $c= 2 \delta^{n^2}$.\footnote{Roughly speaking, one can choose
any $c$, $d$ and $n$ that satisfy the following three properties: 
(1) $d^{1/n}$ is small (equation \eqref{1}),
(2) $d/c$ is small (equation \eqref{11}), 
and (3) $(d/c)^{1/n}$ is not small, say greater than $\delta$.} 
Note that since $ c = 2 \delta^{n^2} < 2 \delta^n \leq 1/2$,
\begin{equation}\label{1}
\frac{4 \delta^{n+1}}{1-c} \leq 8 \delta^{n+1} < 8 \delta^n \leq \frac{1}{2M}.
\end{equation}
Moreover, we also have
\begin{equation}\label{11}
\frac{d}{c-d} = \frac{4 \delta^{n^2+n}}{2 \delta^{n^2} - 4 \delta^{n^2+n}} 
 = \frac{2 \delta^n}{1 - 2 \delta^n} \leq 4 \delta^n < \frac{1}{2M}. 
\end{equation}

Let $\mathbb{D}$ be the unit disc in $\mathbb{C}$, $B:=B(0, d) = \{ z \in \mathbb{C} :
|z| < d \}$, and $A := \mathbb{D} \backslash \overline{B}$. Further, we define
$D := \{ z : (z+c)/(1+cz) \in A \}$, $D_1 := \{ z : z^n \in A  \} = \{ z: d/z^n \in A \}$ and 
$D_2 := \{ z: z^{n^2} \in D \}$; i.e., $D$ is the image of $A$ under the M\"obius
transformation $L (z) : = (z-c)/(1-cz)$, and $D_1$ and $D_2$ are preimages of 
$A$ and $D$ under $h_1 (z) := d/z^n$ and $h_2 (z) := z^{n^2}$, respectively. 
Finally we define the bordered Riemann surface $R$ as following:
\begin{equation}\label{0}
R := \left\{ (z_1, z_2) \in \mathbb{C}^2 : z_1 \in D_1, z_2 \in D_2 \mbox{ and } 
     \frac{z_1^{n} - c}{1 - c z_1^{n}} = z_2^{n^2} \right\}.
\end{equation}
See Figure~1.

\begin{figure}
   \input{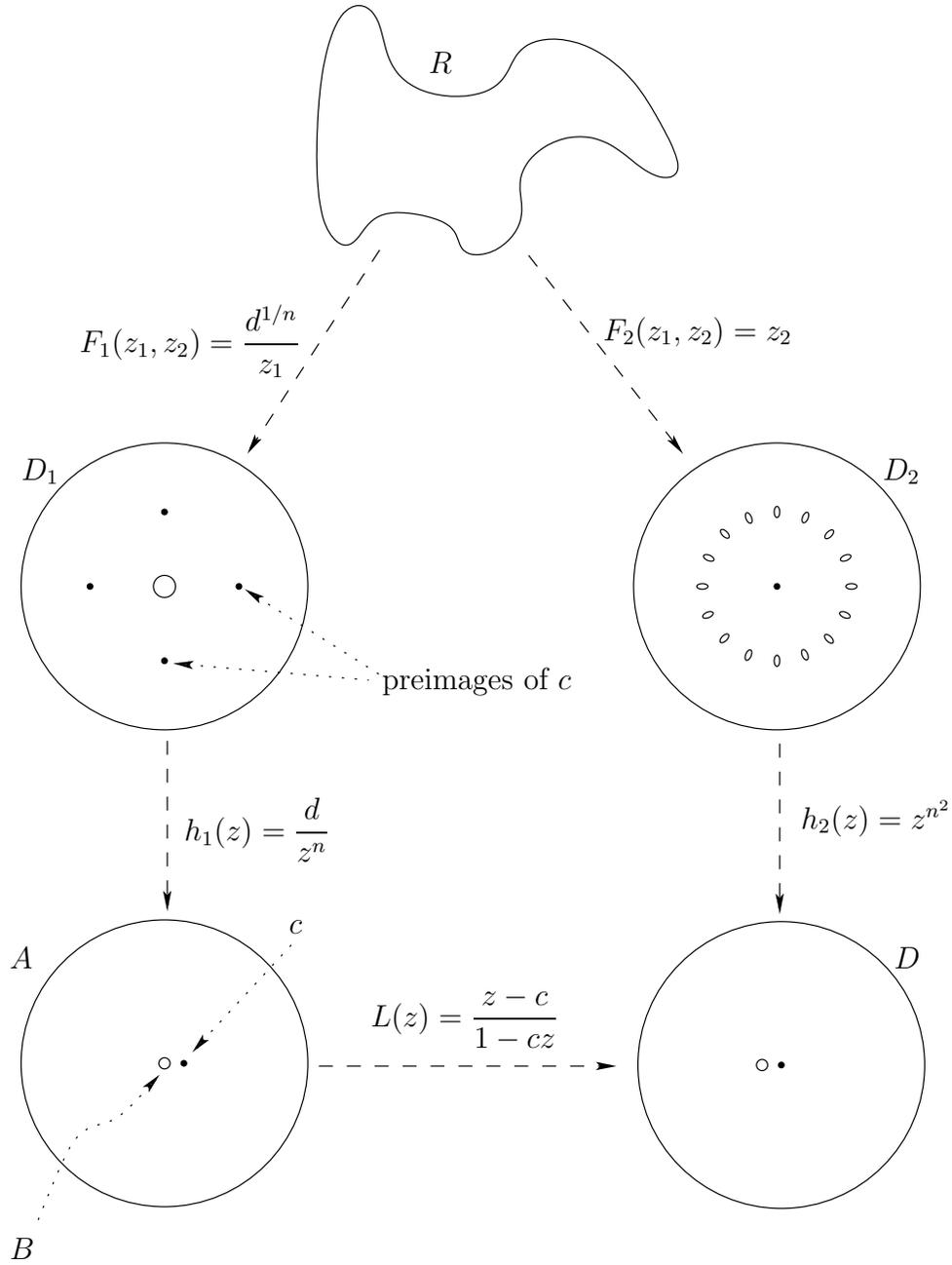}
   \caption{Construction of $R$ for $n=4$}
 \end{figure} 
 
Note that $R$ is an $n$-sheeted covering of $D_2$ and an $n^2$-sheeted branched covering of $D_1$. 
This is because $D_1$ is an $n$-sheeted covering of $A$ and $D_2$ is an $n^2$-sheeted branched covering of $D$.

Now we claim that the Riemann surface $R$, together with the holomorphic functions
\[
F_1 (z_1, z_2) = \frac{d^{1/n}}{z_1} \qquad \mbox{and} \qquad F_2 (z_1, z_2) = z_2,
\]
satisfy the conditions in Theorem~\ref{T}.

First, note that $F_1$ and $F_2$ have values in $D_1$ and $D_2$, respectively. Thus we have
$\max \{ \| F_1 \|, \|  F_2 \| \} \leq 1$.
Furthermore, if $|F_2 (z_1, z_2)| = |z_2| < \delta$, we have
\[
 |z_1^n  - c| = |z_2|^{n^2} |1 - c z_1^n| < 2\delta^{n^2}, 
\]
hence 
\[
|z_1|^n < c + 2 \delta^{n^2} = 4 \delta^{n^2}.
\]
Therefore we have
\[
 |F_1 (z_1, z_2)| = \frac{d^{1/n}}{|z_1|} > \frac{4^{1/n}\delta^{n + 1}}{4^{1/n} \delta^{n}}
                  = \delta,
\] 
and the inequality $\max \{ |F_1 (z_1, z_2)|, |F_2(z_1, z_2)| \} \geq \delta$ holds for all $(z_1, z_2) \in R$; 
i.e., $(F_1, F_2)$ becomes a pair of corona data for the given $\delta$. 

It remains to show that $\max \{ \|G_1\|, \|G_2\| \} \geq M$ for any corona solutions $G_1, G_2$ such that 
\begin{equation}\label{4}
F_1 G_1 + F_2 G_2 = 1.
\end{equation}
In fact, we will show that $\|G_1\| \geq M$. 
To prove this claim, we assume, without loss of generality, that $G_1$ and $G_2$
are holomorphic across the boundary of $R$. Then we define for all $z \in A$,
\[
f(z) := \frac{1}{n^3} \sum F_1(z_1, z_2) G_1(z_1, z_2),
\]
where the summation is over all the points $(z_1, z_2) \in R$ such that $z_1^{n} = z$, counting
multiplicity. (Note that the map $(z_1, z_2) \mapsto z_1^{n}$ is an $n^3$-sheeted branched covering from $R$ to $A$.) Then $f$ is analytic in (a neighborhood of) $A$. 

Since $F_2 (z_1, z_2) = z_2 =0$ when $z_1^n =c$, it is easy to see from \eqref{4} that $f(c) = 1$. On the other hand, 
$|f (z)| \leq  \| G_1 \|$ for all  $z \in A$ since $\|F_1\| \leq 1$, and 
$|f (z)| \leq  4 \delta^{n+1}\| G_1 \|$ for $|z|=1$ since
on $\{ |z_1| = 1 \}$ we have 
\[
|F_1 (z_1, z_2)| = \frac{d^{1/n}}{|z_1|} = 4^{1/n} \delta^{n+1} \leq 4 \delta^{n+1}.
\]
Therefore by Cauchy's integral formula, 
\begin{align*}
1 = |f(c)| & = \left| \frac{1}{2\pi i}\int_{|\xi|=1} \frac{f(\xi)}{\xi - c} d\xi - 
                    \frac{1}{2\pi i} \int_{|\xi|=d} \frac{f(\xi)}{\xi - c } d\xi \right| \\
& \leq \frac{4 \delta^{n+1} \| G_1 \|}{1-c} +  \frac{2 \pi d \| G_1 \|}{2 \pi (c - d)}.
\end{align*}
Now one can easily see that this inequality together with \eqref{1} and \eqref{11} proves the claim.
This completes the proof.

\section{Further remarks} 
\noindent 1. In the construction of $R$, one can take $F_1$ as the projection map $(z_1,z_2) \mapsto z_1$, but
then it is needed to modify the definition \eqref{0} of $R$ to 
\[
R := \left\{ (z_1, z_2) \in \mathbb{C}^2 : z_1 \in D_1, z_2 \in D_2 \mbox{ and } 
     \frac{d/z_1^{n} - c}{1 - c d /z_1^{n}} = z_2^{n^2} \right\}
\]
because we want to make the pair $(F_1, F_2)$  a set of corona data satisfying \eqref{bound}.
\bigskip

\noindent 2. Consider the function $h(z) := z^n$ defined on $D_1$. 
Then it is not difficult to see that 
the Riemann surface $R$ constructed in Section~\ref{proof} is nothing but 
the Riemann surface of the multi-valued function 
\[
h^{-1} \circ L^{-1} \circ h_2 (z) = \left( \frac{z^{n^2} + c}{1 + c z^{n^2}} \right)^{1/n}
\]
defined on $D_2$. Note that it takes values in $D_1$. Similarly, one can consider $R$ as a Riemann surface
of the multi-valued function 
\[
h_2^{-1} \circ L \circ h (z) = \left( \frac{z^{n} - c}{1 - c z^{n}} \right)^{1/n^2}
\]
defined on $D_1$. 
\bigskip

\noindent 3. We can construct $R$ by using cutting and pasting method. For example we can construct the Riemann surface of 
$h^{-1} \circ L^{-1} \circ h_2$ over $D_2$ in the following way: we make $n^2$ cuts
on $D_2$ radially so that each cut connects a hole to the outer boundary of $D_2$ (i.e., to the unit circle). 
We denote this region ($D_2$ minus cuts) by $D(1)$, and enumerate the cuts by $e(1, k, l),$ 
$k=1,\ldots, n^2,$ $l=1,2$ so that $e(1,k,1) = e(1,k,2)$ as sets, and as $z$ approaches to $e(1,k,1)$
the argument of $z$ increases. Let $D(j)$, $j=1, \ldots, n$, be the copies of $D(1)$ with the 
corresponding cuts $e(j,k,l)$, $j=1,\ldots,n$, $k=1,\ldots, n^2, l=1,2$. Now for all $j$ (mod $n$),
we paste $D(j)$ and $D(j+1)$
by identifying $e(j,k,1)$ with $e(j+1, k, 2)$, $k=1,2,\ldots, n^2$. The resulting surface
is conformally equivalent to $R$ with the natural projection map $\pi \approx F_2$. 
Note that by the argument of analytic continuation, 
the map $h^{-1} \circ L^{-1} \circ h_2 \circ \pi$ is well-defined on $R$,
hence analytic. We leave the details to the reader.
\bigskip

\noindent 4. One may end up with the same Riemann surface $R$ from interpolation problems. 
Fix $\epsilon \in (0,1/2)$
and let $D_1' := \{ z : \epsilon < |z| <1 \}$. Choose a natural number $n$ sufficiently large so that
$2^{-n} < \epsilon$, and we let $E_n$ be the set of $n$-th roots of $2^{-n}$.
Note that $|z| = 1/2$ for all $z \in E_n$.

Let us consider the following two interpolation problems: (1) find $G_1 \in H^\infty (D_1')$ 
(with the smallest uniform norm)
such that $G_1 (z) = \overline{z}$ for all $z \in E_n$, and (2) find $F_2 \in H^\infty (D_1')$ 
(with the largest $\delta_0 := \min_{z \in D_1'} \{ |F_2 (z)|, |z| \}$) 
such that $\| F_2 \| =1$ and $F_2 (z) = 0$ for all $z \in E_n$. 

Any solution $G_1$ of (1) has a uniform norm greater than $C/\epsilon$, 
for some absolute constant $C$.
To see this, one may repeat the argument in Section~\ref{proof}; i.e., for $w \in \{ z : \epsilon^n < |z| <1 \}$ define 
\[
f(w) = \frac{1}{n} \sum z G_1 (z),
\]  
where the summation is over all $z \in D_1'$ such that $(\epsilon/z)^n = w$. Note that $f(2^n \epsilon^n) = 1/4$ since 
$z G_1 (z) = 1/4$ for $z \in E_n$, and then Cauchy's integral formula 
gives a lower bound estimate $\| G_1 \| \geq C / \epsilon$.

On the other hand, any solution $F_2$ of (2)
should yield a small $\delta_0 = o(1)$ as $\epsilon \to 0$. To see this, let
$F_1 = z$ and $F_2$ be the solution of (2). Now if $\delta_0 \ne o(1)$, the pair $(F_1, F_2)$
became a set of corona data on $D_1'$ with corresponding 
$0 < \delta \leq \liminf_{\epsilon \to 0} \delta_0$. But then
any corona solutions $G_1$ and $G_2$ such that $F_1 G_1 + F_2 G_2 =1$ 
would have a bound $\geq C/\epsilon$, because
$G_1 / 4$ should be a solution of (1). This violates the corona theorem on
annuli \cite{Sc}, \cite{St}. (In fact, it violates a statement slightly stronger than the corona theorem,
which is true for annuli; i.e., for \emph{any} annulus $D_1'$ and corona data defined
on $D_1'$, there always exist corona solutions with bound $\leq M = M(\delta)$, where $M$ 
does not depend on $D_1'$. See \cite{Gam}, p.~47, for details.) Therefore to make
$F_1$ and $F_2$ corona data, or to get a solution for (2) with large $\delta_0$, 
we take a number $N$ such that the multi-valued function,
\[
F (z) = \left( \frac{z^n-2^{-n}}{1-2^{-n} z^n} \right)^{1/N},
\]
has a modulus $\geq 1/4$ for $|z| < 1/4$. (Such $N$ should be approximately greater than $n^2$ as in Section~\ref{proof}. Also note that $F^N$ is a solution for (2).) 
Now since $F$ is not analytic
on $D_1'$, we consider the Riemann surface of $F$ over $D_1'$, which gives us the Riemann surface $R$
constructed in Section~2 (with $\delta = 1/4$).

\end{document}